%% file: main.tex
\documentclass[graybox]{SNmult}
\usepackage{type1cm}
\usepackage{makeidx}
\usepackage{graphicx}
\usepackage{multicol}
\usepackage{multirow}
\usepackage{mathtools}
\usepackage[bottom]{footmisc}
\usepackage{amsmath,amssymb}
\usepackage{tabularx}
\usepackage{booktabs}
\usepackage{makecell}
\usepackage{microtype}
\usepackage{placeins}

\usepackage{newtxtext}
\usepackage[varvw]{newtxmath}
\usepackage{comment}
\usepackage[pdftex,
            colorlinks=true,
            linkcolor=black,
            citecolor=green,
            filecolor=black,
            urlcolor=green]{hyperref}
\usepackage{caption}
\makeindex

\begin{document}
\title*{A Parallel Heuristic for the Fixed Cardinality Binary Quadratic Program}
\titlerunning{A parallel heuristic for FCBQP}
\author{Thi Thanh Tu Le\orcidID{0009-0001-8161-4535} and\\
Thorsten Koch\orcidID{0000-0002-1967-0077}}

\institute{
Thi Thanh Tu Le \at Software and Algorithms for Discrete Optimization, Technische Universit\"at Berlin, Germany; and
Applied Optimization Department, Zuse Institute Berlin, Germany;
\email{thi.le@zib.de}
\and 
Thorsten Koch \at
 Chair of Software and Algorithms for Discrete Optimization, Technische Universit\"at Berlin, Germany; and
Applied Optimization Department, Zuse Institute Berlin, Germany;
\email{koch@zib.de}
}
\maketitle
\abstract{The Fixed Cardinality Binary Quadratic Program (FCBQP) is an optimization problem that requires a prescribed number of variables to take the value 1.
We propose a parallel CPU-based heuristic that combines impact-based variable
selection, tabu search, guided restarts, and a shared elite pool. Computational
experiments on 342 benchmark instances from the maximum dispersion problem,
$k$-clique problem, and minimum graph bisection problem show that the new algorithm
is competitive and can reach best-known or
reference values on instances with more than 20000 variables within one minute.}

\section{Introduction}
\label{sec:introduction}
We consider the \textit{Fixed Cardinality Binary Quadratic Program} (FCBQP)
\begin{equation}
\label{eq:fcbqp}
\min \; f(x)=x^\top Qx
\quad \text{subject to}\quad
\sum_{i=1}^{n}x_i=k,\quad x\in\{0,1\}^{n},
\end{equation}
where $Q=(q_{ij})_{1 \leq i , j \leq n} \in\mathbb{R}^{n\times n}$ is the cost matrix, $k \in \mathbb{N}$, and
$0 \leq k\leq n$. The cardinality constraint enforces that exactly $k$ variables take the value 1 in each feasible solution. W.l.o.g., we assume that $Q$ is symmetric, i.e., $q_{ij} = q_{ji}$ for all $1 \leq i < j \leq n$. Moreover, by complementing variables and transforming the cost matrix, we assume $0 \leq k \leq \lfloor \frac{n}{2} \rfloor$.

The FCBQP is NP-hard and generalizes several combinatorial optimization problems, including the maximum dispersion problem (MDP)~\cite{zhou2017opposition} and graph bisection problem (GBP)~\cite{walshaw2002graph}. Furthermore, any unconstrained binary quadratic program can be solved by considering FCBQPs for all cardinalities $k = 0, \ldots, n$, with the $k$-clique problem as a special case~\cite{sotirov2020solving}. Other applications of the FCBQP include facility location~\cite{zhou2017opposition}, molecular
structure design~\cite{zhou2017opposition}, clustering~\cite{sotirov2020solving}, and wind farm layout optimization \cite{cazzaro2022variable}.

Many heuristics have been proposed to efficiently find high-quality solutions for the
FCBQP~\cite{zhou2017opposition,
hochbaum2023breakpoints, liu2022two, nogueira2025gpu}, but most of them focus on the MDP. We propose a new heuristic for large and challenging FCBQP instances, together with a parallel version that improves robustness and solution quality. The algorithmic details and computational results are presented in Sections~\ref{sec:algorithm} and~\ref{sec:computation}, respectively, followed by the conclusions in Section~\ref{sec:conclusion}.

\section{Proposed Heuristic for the FCBQP} \label{sec:algorithm}
A feasible solution of~\eqref{eq:fcbqp} contains $k$ \textit{selected} variables ($x_i=1$) and $n-k$ \textit{unselected} variables ($x_i=0$). The trivial cases $k = 0, k =1$, $k = n-1$, and $k =n$ are handled separately. In all other cases, the algorithm first constructs an initial candidate set using \textit{impact-based variable selection}. It then generates feasible solutions by \textit{random-greedy construction}, \textit{random-diverse restart}, \textit{elite-based restart}, and \textit{path-relinking}. Each solution is then improved by \textit{tabu search}, while an \textit{elite pool} stores high-quality and diverse solutions. 

\textit{2.1. Impact-based variable selection.}
For each variable $i$, the interaction impact score is defined by
$\operatorname{imp}(i) = q_{ii} + 2 \frac{k-1}{n-1}\sum_{j\ne i}q_{ij}
$. This represents the expected contribution of variable $i$ if $i$ is selected and the remaining $k-1$ selected variables are sampled uniformly from the other $n-1$ variables. An initial set $C_1$ of potential candidates is built by selecting $v_0$ variables with the lowest interaction impact scores.

\textit{2.2. Tabu search.}
At each iteration, a swap exchanges one selected variable with one unselected variable. 
After a variable’s value changes, it becomes tabu for several iterations, unless the move satisfies the \textit{aspiration} criterion, i.e., improves the best solution found during the current tabu search. 
The tabu tenure, denoted by $T$, is adaptive and is set as shown in Table \ref{tab:parameters} based on the number of consecutive non-improving iterations $T_{\operatorname{nimp}}$ and a factor $t_\alpha$.  After a swap $(i,j)$, where
$x_i : 1 \to 0$ and $x_j : 0 \to 1$, variable $i$ is  tabu for $T$ iterations, while variable $j$ is tabu for
$\lfloor 2T/3 \rfloor$ iterations. A swap is called \textit{admissible} if neither of the two variables is tabu or if the aspiration criterion holds.
To avoid evaluating all possible swaps, we use the support value
$h(i)=\sum_{j=1}^{n}q_{ij}x_j, \text{for } i=1,\ldots,n.$
Only the $v_1$ selected variables with the largest $h(i)$ and the $v_2$ unselected variables with the smallest $h(i)$ are considered for swaps. For a candidate swap $(i,j)$ with $x_i:1 \to 0$ and $x_j:0 \to 1$, the objective change is $\Delta(i,j)=2(h(j)-h(i))+q_{ii}+q_{jj}-2q_{ij}$.
The best admissible swap is applied; if none exists, the best candidate swap is chosen ignoring tabu status. The support values are updated as $h(l) \coloneqq h(l) -q_{li} + q_{lj}$ for $l = 1, \dots, n$. The objective value and tabu list are also updated after each swap. The tabu search terminates after $T_{\max}$ iterations or after 2000 consecutive non-improving iterations.

\textit{2.3. Elite pool and potential candidate set management.}
The algorithm maintains a global elite pool of at most 30 solutions. A new solution is accepted only if it is distinct from the solutions in the elite pool. If the pool is not full, the new solution is inserted when its minimum Hamming distance to the pool is at least $d$, or at least $d'=\max(2,\lceil0.8d\rceil)$ if it improves the best elite solution. 

When the elite pool is full, the algorithm may replace an elite solution
$e_r$. Let $\delta_r$ be the minimum Hamming distance between the new solution
and the remaining elite solutions. A replacement is admissible if the new
solution either (i) is no worse than the worst elite solution and
$\delta_r \geq d$, (ii) has a better objective value than $e_r$, satisfies
$\delta_r \geq d'$, and preserves at least $80\%$ of the diversity of $e_r$,
or (iii) improves the best elite solution and satisfies $\delta_r \geq d'$.
Among all admissible replacements, the one giving the largest objective
improvement is selected, with greater diversity used as a secondary criterion.

The algorithm also constructs a candidate set $C_2$.
This set consists of all variables in $C_1$ and all variables that take the value 1 in at least one solution in the elite pool.

\textit{2.4. Guided restarts.}
The algorithm uses the following four construction strategies.

\textit{Random-greedy construction.} This strategy first randomly selects $v_3$ variables from $C_2$ and fixes them to 1. The remaining selected variables are added iteratively from among all currently unfixed variables, using the marginal contribution $q_{ii}+2h(i)$; at each step, one variable is chosen at random from the 20 candidates with the smallest marginal contributions. If the elite pool is nonempty, the construction is repeated until the generated solution has minimum Hamming distance of at least $d$ from the elite pool, or until 200 attempts are reached, in which case the most distant solution is returned.

\textit{Random-diverse restart.} 
Up to 200 feasible solutions are generated randomly, and the first solution
whose minimum Hamming distance from the elite pool is at least $d$ is returned. If none
satisfies this condition, the most distant generated solution is used.

\textit{Elite-based restart.} Variables that are selected in at least 85\% of the solutions in the current elite pool are candidates for fixing to 1. At most $\lfloor0.85k\rfloor$ such variables are fixed, prioritizing higher elite frequencies and breaking ties by smaller variable indices. The remaining variables are selected randomly,
again using at most 200 attempts to satisfy the distance threshold to the elite pool.

\textit{Path-relinking.} 
A path between two feasible solutions is a sequence of swaps that gradually
transforms one solution into the other, where at each step, one variable
selected only in the current solution is replaced by one selected only in the target solution. The current incumbent serves as the initiating solution. 
 The guiding solution is selected by balancing objective quality and Hamming distance. Whenever possible, its Hamming distance from the incumbent must be at least $d$; otherwise, the most distant elite solution is selected. Along each path, restricted-candidate swaps are constructed using $h(i)$, with each move selected uniformly at random from among the 10 best swaps. The procedure is applied in both directions, and the best intermediate solution is returned.

\textit{2.5. Parallelization and restart ordering.} 
The algorithm uses the C++ Parallel STL to run multiple search lanes concurrently. Each lane has its own random seed and incumbent, and, by default, all lanes share a global elite pool. At each restart, a lane takes a snapshot of the elite pool. If the elite pool is nonempty, path-relinking is used after 10 consecutive non-improving starts; otherwise, the restart strategy is chosen according to the restart index: an elite-based restart if an incumbent exists and the index is divisible by 10, a random-diverse restart if it is congruent to 9 modulo 10, and random-greedy construction otherwise. The generated solution is then improved by tabu search and considered for insertion into the elite pool. At the end of each epoch, candidate solutions are merged into the pool in lane-ID order. The search is deterministic for fixed parameters when no runtime limit is
used. For these experiments, we enable early stopping when the best-known or reference value is reached.
The parameter settings are in Table~\ref{tab:parameters}.
\begin{table}[h]
\centering
\caption{Parameter settings.}
\label{tab:parameters}
\scriptsize
\renewcommand{\arraystretch}{1.25}
\setlength{\tabcolsep}{3pt}

\begin{tabular}{|p{0.35\textwidth}|p{0.6\textwidth}|}
\hline

\begin{minipage}[t]{\linewidth}
\vspace{0pt}
\begin{tabularx}{\linewidth}{@{}llX@{}}
\textbf{Par.} & \textbf{Value}\\
\hline
$v_0$ 
& $\min\{k,\lceil0.4n \rceil\}$\\

$d$ 
& $\lceil 0.2k \rceil$ \\

$T$ 
& $\begin{cases}
t_\alpha, & T_{\mathrm{nimp}} = 0,\\
2t_\alpha, & 0 <T_{\mathrm{nimp}} < 100,\\
3t_\alpha, & \text{otherwise}.
\end{cases}$ \\

$T_{\max}$ 
& $\begin{cases}
10000, & n \le 7000,\\
20000, & \text{otherwise}.
\end{cases}$ \\

$v_3$ & $\min\{|C_2|, \lfloor 0.67k + 0.5 \rfloor\}$ 

\end{tabularx}
\end{minipage}
&
\begin{minipage}[t]{\linewidth}
\vspace{0pt}
\begin{tabularx}{\linewidth}{@{}llX@{}}
\textbf{Par.}  & \textbf{Value} \\
\hline
$v_1$ 
& $\begin{cases}
\min\{\max\{80,\lfloor k/3\rfloor\},100, k\}, & n \le 7000,\\
\min\{\max\{80,\lfloor k/3\rfloor\},300, k\}, & \text{otherwise}.
\end{cases}$ \\

$v_2$ 
& $\begin{cases}
\min\{\max\{120,\lfloor(n-k)/12 \rfloor\},200, n-k\}, & n \le 7000,\\
\min\{\max\{120,\lfloor (n-k)/12\rfloor \},300, n-k\}, & \text{otherwise}.
\end{cases}$ \\

$t_\alpha$ 
& $\begin{cases}
3, & n \le 500,\\
5, & 500 <n \le 2000,\\
7, & 2000 < n \le 5000,\\
9, & 5000 <n \le 8000,\\
15, & 8000 <n \le 20000,\\
\max\{1,\lfloor 0.2\sqrt{\min\{k,n-k\}} \rfloor\}, & \text{otherwise}.
\end{cases}$ \\
\end{tabularx}
\end{minipage}
\\
\hline
\end{tabular}
\end{table}

\FloatBarrier
\vspace*{-15pt}
\section{Computational Results} \label{sec:computation}
As mentioned earlier, the FCBQP encompasses several combinatorial optimization problems, including MDP, $k$-clique, and GBP. We test our algorithm on benchmark instances from these subproblems, performing 5 runs per instance with a 60-second time limit per run.
All experiments were run on a machine equipped with a 24-core AMD Ryzen Threadripper PRO 7965WX processor and an NVIDIA GeForce RTX 4090 GPU. The source code, detailed computational results, and datasets are available in GitHub repositories%
\footnote{\url{https://github.com/LeThiThanhTu/FCBQP_parallel}}%
\footnote{\url{https://github.com/LeThiThanhTu/Benchmark_instances_FCBQP}}.
The code was compiled
with GCC 12 in C++20 mode, with link-time
optimization and native-architecture optimization enabled. The GNU
Parallel STL used the TBB backend. The executable was built using
\texttt{make build-gcc-fast}.
 We report the best gap (Best gap), average gap (Avg. gap), number of 
runs reaching the best-known or reference value (\#Succ), and average time-to-solution (Avg.\ TTS). For each  instance in the minimization form \eqref{eq:fcbqp}, the average gap is the mean relative deviation
from the best-known or reference value over 5 runs, i.e.,
$\mathrm{Avg. \ gap}
=
\frac{1}{5}
\sum_{r=1}^{5}
\frac{ f_r-f^\ast}{\max\{1,\lvert f^* \rvert\}}$,
where \(f_r\) is the best objective value obtained in the $r$-th run,
and $f^\ast$ is the best-known or reference value. For each dataset, Avg. gap is the average of
these values over all instances in the dataset. 
Best gap is the average, over all instances in the dataset,
of the relative deviation between the best objective value found over 5 runs and
the best-known or reference value. 
For each instance, Avg. TTS is averaged over successful
runs only and is reported as ``--'' when there is no successful run.
Bold values indicate the best results.

The MDP aims to select exactly $k$ elements from a given set so that the total pairwise distance among the selected elements is maximized. It can be formulated as
$\max \left\{ \frac{1}{2} x^\top D x \;\middle|\; \sum_{i=1}^n x_i = k,\; x \in \{0,1\}^n \right\}$,
where $D=(d_{ij})$ is the pairwise distance matrix and $d_{ii}=0$ for all $i=1,\dots,n$. To convert this to the minimization form
\eqref{eq:fcbqp}, we set $Q=-\frac{1}{2}D$. %and report the original maximization objective value $\frac{1}{2}x^\top Dx$. This corresponds to using the \texttt{max} objective setting with the \texttt{half} output scale in our code.
Our computational assessment uses 160 MDPLIB instances~\cite{ liu2022two, marti2021mdplib} and 90 sparse instances from~\cite{hochbaum2023breakpoints}, for which we impose cardinality constraints $k=100,200,300,400,500$. The best-known values for the MDPLIB instances are taken from the literature~\cite{liu2022two, marti2021mdplib}.  For the instances from~\cite{hochbaum2023breakpoints}, no best-known values are available. Therefore, 
for each instance from~\cite{hochbaum2023breakpoints}, we use as the reference value the best among the values obtained by Gurobi \cite{Gurobi} in 1800 seconds, by OBMA~\cite{zhou2017opposition} and our algorithm in 600 seconds, and by the breakpoint-based method in~\cite{hochbaum2023breakpoints}.
We compare our algorithm with GPU-TS-f3, a GPU-based tabu search algorithm considered one of the best heuristics for the MDP~\cite{nogueira2025gpu}. For a fair comparison, we modified the \texttt{GPU-TS-f3.cu} code\footnote{\url{https://github.com/WilliamRosendo/GPU-Tabu-Search}} so
that, like our algorithm, it stops upon reaching the best-known or
reference value. Its stopping test uses the same numerical tolerance as our algorithm, thereby avoiding delayed stopping caused by floating-point
rounding. We ran GPU-TS-f3 and our algorithm on each instance 5 times, with a 60-second limit per run, and recorded the best solution from each run. The MDP test sets are divided into Group I, in which all methods reach the best-known or reference value in every run, and the more difficult Group II. The absolute objective tolerance is set to $10^{-5}$.

Group I contains five MDP datasets (20 instances each): SOM-b, GKD-c, MDG-a\_n500, MDG-b\_n500, and LA. Both the 1-core and 16-core versions of our algorithm, as well as GPU-TS-f3, reach the best-known values for all instances in all 1500 runs, with an average time-to-solution of 0.21 seconds. We omit a  comparison of these methods for Group I because their TTS are negligible. Results of Group II are reported in Table~\ref{tab:mdp_2}.

Table~\ref{tab:mdp_2} shows that our 16-core algorithm provides the best overall performance on these challenging MDP datasets. GPU-TS-f3 performs better than our 1-core algorithm on MDG-a\_n2000, MDG-b\_n2000, b2500, MDG-c, p3000-5000, and  StackOverflow, but its performance deteriorates substantially on IMDB, Bibsonomy, SyntheticTF, and DBLP. Note that the breakpoint-based method\footnote{Since MDP is a special case of QKP, we ran the breakpoint-based QKP implementation at 
\url{https://github.com/phil85/breakpoints-algorithm-for-qkp}
to compute the reported gaps.} in~\cite{hochbaum2023breakpoints, hochbaum2025fast} runs within a few seconds but yields average gaps of 0.02 on IMDB, 0.015 on StackOverflow, 0.02 on Bibsonomy, 0.04 on SyntheticTF, and 0.036 on DBLP, all larger than those of our serial version. 
\begin{table}[h]
\centering
\caption{Results on the MDP datasets in Group II.}
\label{tab:mdp_2}
\scriptsize
\setlength{\tabcolsep}{3pt}
\renewcommand{\arraystretch}{0.85}

\resizebox{\textwidth}{!}{%
\begin{tabular}{lccccccccc}
\toprule
&
\multicolumn{3}{c}{New alg. (16 cores)}
&
\multicolumn{3}{c}{New alg. (1 core)}
&
\multicolumn{3}{c}{GPU-TS-f3} \\
\cmidrule(lr){2-4}
\cmidrule(lr){5-7}
\cmidrule(lr){8-10}
Dataset
& Best gap
& Avg. gap
& \#Succ
& Best gap
& Avg. gap
& \#Succ
& Best gap
& Avg. gap
& \#Succ \\
\midrule

MDG-a\_n2000
& \textbf{0} & \textbf{$<10^{-5}$} & \textbf{99/100}
& 0 & 0.00007 & 65/100
& 0 & $<10^{-5}$ & 88/100 \\

MDG-b\_n2000
& \textbf{0} & \textbf{0} & \textbf{100/100}
& $<10^{-5}$ & 0.00007 & 50/100
& 0 & $<10^{-5}$ & 79/100 \\

b2500
& \textbf{0.00001} & \textbf{0.00001} & \textbf{45/50}
& 0.00005 & 0.00023 & 9/50
& 0.00001 & 0.00004 & 34/50 \\

MDG-c
& \textbf{0.00019} & \textbf{0.00021} & \textbf{37/100}
& 0.00029 & 0.00039 & 3/100
& 0.00021 & 0.00026 & 13/100 \\

p3000-5000
& \textbf{$<10^{-5}$} & \textbf{$<10^{-5}$} & \textbf{22/50}
& 0.00001 & 0.00003 & 4/50
& $<10^{-5}$ & 0.00001 & 16/50 \\

StackOverflow
& \textbf{0} & \textbf{0} & \textbf{25/25}
& 0 & 0.00117 & 16/25
& 0 & 0.00021 & 19/25 \\

IMDB
& \textbf{0} & \textbf{$<10^{-5}$} & \textbf{24/25}
& 0.00008 & 0.00028 & 14/25
& 0.00244 & 0.00632 & 1/25 \\

Bibsonomy
& \textbf{0.00056} &  \textbf{0.00076} & \textbf{9/25}
& 0.00156 & 0.00263 & 6/25
& 0.15776 & 0.22716 & 0/25 \\

SyntheticTF
& \textbf{0.00057} & \textbf{0.00097} & \textbf{77/250}
&  0.00270 & 0.00464&   27/250
& 0.25266 &  0.28796 & 0/250 \\

DBLP
& \textbf{0.00408} & \textbf{0.00473} & 0/25
& 0.00823 &  0.01166 & 0/25
& 0.31018 & 0.34629 & 0/25 \\

\bottomrule
\end{tabular}%
}
\end{table}

Given an undirected graph $G=(V,E)$, the $k$-clique problem asks if there are $k$ nodes in the graph that form a clique. For any fixed size $k$, $G$ contains a clique of size $k$ iff
$\max \left\{x^\top Ax \;\middle|\;
\sum_{i=1}^{n}x_i=k,\ x\in\{0,1\}^{n}\right\}=k(k-1),$ where $A$ is the adjacency matrix of $G$. To express this in the minimization form \eqref{eq:fcbqp}, we set $Q = -A$.  
For our computational experiments, we use three benchmark test sets originally generated for the Maximum Clique Problem (MCP): 7 instances used in ~\cite{sotirov2020solving} based on real-world graphs, 41 BHOSLIB instances %\footnote{\url{https://rb-bench.github.io/benchmarks/hidden-opt-graph.html}}  
\cite{bhoslib}, and 37 DIMACS instances
%\footnote{\url{https://archive.dimacs.rutgers.edu/pub/challenge/graph/benchmarks/volume/Clique/}}
\cite{dimacs}. 
%For each instance, we consider the corresponding $k$-clique problem, where $k$ is either the best-known maximum clique size reported in the literature or the maximum clique size found by Gurobi~\cite{Gurobi}.
Table~\ref{tab:mcp_test_1} reports the results on the MCP instances used in~\cite{sotirov2020solving}. 
We remove isolated vertices from each graph and use Gurobi to determine
the maximum clique size of the resulting graph.
Both the 1-core and 16-core versions of our algorithm confirm the clique sizes for every instance in all runs.
The 16-core version is particularly effective on the largest MCP instances, such as ca-AstroPh  and ca-CondMat, for which  it confirms the clique sizes
in approximately 2.3 seconds on average. Note that the clique sizes  here are much larger than those reported in~\cite{sotirov2020solving},
and our algorithm confirms these larger clique sizes much
faster than the random coordinate descent method
in~\cite{sotirov2020solving}. 

\begin{table}[h]
\centering
\scriptsize
\setlength{\tabcolsep}{3pt}
\renewcommand{\arraystretch}{1.00}
\captionsetup{skip=2pt}
\caption{Results on the MCP instances used in~\cite{sotirov2020solving}. Here, $n$ denotes the number of vertices after removing isolated vertices,
and ``Best obj.'' reports the original  objective value 
$x^\top A x=k(k-1)$.}
\label{tab:mcp_test_1}

\begin{tabular*}{\textwidth}{
    @{\extracolsep{\fill}}
    l r r r cc cc cc
}
\toprule
& & & &
\multicolumn{2}{c}{New alg. (16 cores)}
& \multicolumn{2}{c}{New alg. (1 core)}
& \multicolumn{2}{c}{GPU-TS-f3} \\
\cmidrule(lr){5-6}
\cmidrule(lr){7-8}
\cmidrule(lr){9-10}
Instance & $n$ & $k$ & Best obj.
& Avg. TTS & \#Succ
& Avg. TTS & \#Succ
& Avg. TTS & \#Succ \\
\midrule
jazz       & 198   & 30  & 870
& 0.0005 & 5/5 & \textbf{0.0002} & \textbf{5/5} & 0.0019 & 5/5 \\
email      & 1133  & 12  & 132
& \textbf{0.0019} & \textbf{5/5} & 0.0082 & 5/5 & 0.0873 & 5/5 \\
ca-GrQc    & 5241  & 44  & 1892
& 0.0006 & 5/5 & \textbf{0.0003} & \textbf{5/5} & 0.0034 & 5/5 \\
ca-HepTh   & 9875  & 32  & 992
& \textbf{0.6145} & \textbf{5/5} & 3.6829 & 5/5 & 37.0194 & 4/5 \\
ca-HepPh   & 12006 & 239 & 56882
& 0.0044 & 5/5 & \textbf{0.0033} & \textbf{5/5} & 0.0161 & 5/5 \\
ca-AstroPh & 18772 & 57  & 3192
& \textbf{2.2037} & \textbf{5/5} & 10.3454 & 5/5 & 23.6059 & 3/5 \\
ca-CondMat & 23133 & 26  & 650
& \textbf{2.3492} & \textbf{5/5} & 15.8320 & 5/5 & 17.8467 & 5/5 \\
\bottomrule
\end{tabular*}
\end{table}

Table~\ref{tab:mcp_test_2_3} presents the DIMACS and BHOSLIB results.
Tables~\ref{tab:mcp_test_1} and~\ref{tab:mcp_test_2_3} show that our
1-core version outperforms GPU-TS-f3 on all of these MCP test sets, while the 16-core version further
improves wall-clock TTS and success rates on many instances.

\begin{table}[h]
\caption{Results on the BHOSLIB and DIMACS MCP datasets.}
\label{tab:mcp_test_2_3}
\scriptsize
\setlength{\tabcolsep}{3pt}
\begin{tabular*}{\textwidth}{
    @{\extracolsep{\fill}}
    l ccc ccc ccc
}
\toprule
&
\multicolumn{3}{c}{New alg. (16 cores)}
& \multicolumn{3}{c}{New alg. (1 core)}
& \multicolumn{3}{c}{GPU-TS-f3} \\
\cmidrule(lr){2-4}
\cmidrule(lr){5-7}
\cmidrule(lr){8-10}
Dataset
& Best gap & Avg. gap & \#Succ
& Best gap & Avg. gap & \#Succ
& Best gap & Avg. gap & \#Succ \\
\midrule
BHOSLIB
& \textbf{0.0002} & \textbf{0.0004} & \textbf{111/205}
& 0.0005 & 0.0009 &  68/205
& 0.0012 & 0.0016 & 33/205 \\

DIMACS
& \textbf{0.0022} & \textbf{0.0024} & \textbf{142/185}
& 0.0024 & 0.0026 & 135/185
& 0.0025 & 0.0028 & 127/185 \\
\bottomrule
\end{tabular*}
\end{table}

The GBP seeks to partition a graph $G=(V,E)$ into two equal or nearly equal parts while minimizing the number of edges between them. It can be formulated as
$\min\left\{
\sum_{(i,j)\in E}\left(x_i+x_j-2x_i x_j\right)
\;\middle|\;
\sum_{i=1}^{n}x_i=k,\; x\in\{0,1\}^{n}
\right\}$, where $k = \lfloor \frac{n}{2} \rfloor$ and $x_i=1$ if vertex $i$ is assigned to one part and $x_i=0$ otherwise. 
%For an even number of vertices, $k=n/2$; for an odd number of vertices, one may take $k=\lfloor n/2\rfloor$, which yields two parts of sizes $\lfloor n/2\rfloor$ and $\lceil n/2\rceil$. 
Hence, the GBP can be written in the
form \eqref{eq:fcbqp} by setting
$q_{ii}=\deg(i)$, $q_{ij}=q_{ji}=-1$ for each edge $(i,j)\in E$, and
$q_{ij}=q_{ji}=0$ for $(i,j)\notin E$.
We use challenging instances from the Walshaw graph partitioning
repository~\cite{walshaw2002graph} to evaluate the performance of our method
on GBP instances.

\noindent
\begin{minipage}[t]{0.53\textwidth}
\vspace{0pt}
Table~\ref{tab:graph_partition} reports the results for several instances
with more than 10000 variables. The 16-core version outperforms
the serial version on these instances. In particular,  the 1-core version fails to reach the best-known value on instances
fe\_4elt2.graph ($n =11143$), fe\_sphere.graph ($n=16386$), and
bcsstk30.graph ($n = 28924$) in all runs,  whereas the 16-core version succeeds in several runs.
\end{minipage}
\hfill
\begin{minipage}[t]{0.44\textwidth}
\vspace{0pt}
\centering
\captionsetup{font=scriptsize,skip=3pt}
\captionof{table}{Results on GBP instances.}
\label{tab:graph_partition}

\scriptsize
\setlength{\tabcolsep}{1.5pt}
\renewcommand{\arraystretch}{1.00}

\begin{tabular}{lcccc}
\toprule
& \multicolumn{2}{c}{New alg. (16 cores)}
& \multicolumn{2}{c}{New alg. (1 core)} \\
\cmidrule(lr){2-3}
\cmidrule(lr){4-5}
Instance
& Avg. TTS & \#Succ
& Avg. TTS & \#Succ \\
\midrule
crack.graph
& \textbf{7}  & \textbf{4/5}
& 21          & 1/5 \\

fe\_4elt2.graph
& \textbf{20} & \textbf{5/5}
& --          & 0/5 \\

bcsstk29.graph
& \textbf{1}  & \textbf{5/5}
& 14          & 5/5 \\

4elt.graph
& \textbf{6}  & \textbf{5/5}
& 9           & 1/5 \\

fe\_sphere.graph
& \textbf{35} & \textbf{1/5}
& --          & 0/5 \\

cti.graph
& \textbf{4}  & \textbf{5/5}
& 42          & 1/5 \\

bcsstk30.graph
& \textbf{6}  & \textbf{2/5}
& --          & 0/5 \\
\bottomrule
\end{tabular}
\end{minipage}

\section{Conclusion}
\label{sec:conclusion}

\vspace{-0.3cm}
We proposed a parallel CPU-based heuristic for the FCBQP that combines impact-based variable selection, tabu search, guided restarts, and a shared elite pool. Computational experiments on benchmark instances from the maximum dispersion, $k$-clique, and graph bisection problems show that the
algorithm is competitive. Under the same wall-clock time limit, the 16-core version
generally achieves higher success rates and better solution quality than the 1-core version, particularly on challenging instances. 
More specifically, the overall success rate increases from 54.97\% (940/1710 runs) with 1 core to 73.27\% (1253/1710 runs) with 16 cores, representing an improvement of 18.30 percentage points. Furthermore, for the 162 instances on which both versions reach the best-known or reference value in all 5 runs, the 16-core version achieves a geometric mean speedup of 7.25 times over the 1-core version in terms of time to solution.

\begin{acknowledgement}
This work is supported by the German Academic Exchange Service (DAAD). The computations were performed on a computer at the Zuse Institute Berlin. 
OpenAI Codex was used to assist with code implementation. All algorithmic design decisions were made by the authors.
\end{acknowledgement}
\vspace{-20pt}
\input{references}

\end{document}

%% file: references.tex
%%%%%%%%%%%%%%%%%%%%%%%% referenc.tex %%%%%%%%%%%%%%%%%%%%%%%%%%%%%%
% sample references
% %
% Use this file as a template for your own input.
%
%%%%%%%%%%%%%%%%%%%%%%%% Springer-Verlag %%%%%%%%%%%%%%%%%%%%%%%%%%
%
% BibTeX users please use
% \bibliographystyle{}
% \bibliography{}
%
%\biblstarthook{